\documentclass[11pt]{amsart}

\usepackage[english]{babel}

\usepackage{amsmath, amsthm, amssymb, amsfonts, amsopn}
\usepackage[foot]{amsaddr}

\usepackage{mathrsfs, mathtools, bm}

\usepackage{array, multirow, makecell, booktabs}
\setcellgapes{3pt}
\renewcommand\arraystretch{1.3}

\usepackage[ruled,vlined]{algorithm2e}

\usepackage{graphicx}

\usepackage[skip=2pt,font=scriptsize]{caption}
\usepackage{subcaption}

\usepackage{enumerate}
\usepackage{enumitem}

\usepackage{xcolor}

\usepackage{hyperref}
\usepackage{cleveref}
\hypersetup{
    colorlinks,
    linkcolor={red!80!black},
    citecolor={green!50!black}
}

\usepackage[margin=1in]{geometry}

\theoremstyle{plain}

\theoremstyle{definition}

\numberwithin{equation}{section}
\numberwithin{figure}{section}
\numberwithin{table}{section}

\newcommand{\R}{\mathbb{R}}

\newcommand{\D}{\mathcal{D}}

\renewcommand{\Sigma}{\boldsymbol{\sigma}}
\renewcommand{\L}{\mathsf{L}}

\renewcommand{\u}{\boldsymbol{u}}

\title[HYCO: Hybrid-Cooperative Learning]{HYCO: A FORMALISM FOR
HYBRID-COOPERATIVE PDE MODELLING}

\author[L. Liverani, E. Zuazua]{Lorenzo Liverani$^\dagger$}
\address{$^\dagger$Chair for Dynamics, Control, Machine Learning, and Numerics (Alexander von Humboldt Professorship), Department of Mathematics, Friedrich-Alexander-Universit\"at Erlangen-N\"urnberg, 91058 Erlangen, Germany.}

\author{Enrique Zuazua$^\dagger$$^*$$^\S$}
\address{$^*$Departamento de Matem\'aticas, Universidad Aut\'onoma de Madrid, 28049 Madrid, Spain.}
\address{$^\S$Chair of Computational Mathematics, Universidad de Deusto. Av. de las Universidades, 24, 48007 Bilbao, Basque Country, Spain.}
\email{lorenzo.liverani@fau.de}
\email{enrique.zuazua@fau.de}

\begin{document}

\begin{abstract}
We present Hybrid-Cooperative Learning (HYCO), a hybrid modeling framework that integrates physics-based and data-driven models through mutual regularization. Unlike traditional approaches that impose physical constraints directly on synthetic models, HYCO treats both components as co-trained agents nudged toward agreement. This cooperative scheme is naturally parallelizable and demonstrates robustness to sparse and noisy data. Numerical experiments on static and time-dependent benchmark problems show that HYCO can recover accurate solutions and model parameters under ill-posed conditions. The framework admits a game-theoretic interpretation as a Nash equilibrium problem, enabling alternating optimization. This paper is based on the extended preprint \cite{XXX}.
\end{abstract}


\maketitle

\section{Introduction}
\label{sec:intro}

\noindent
Mathematical modeling has long relied on two dominant paradigms. Physics-based models are grounded in physical laws and typically involve a small number of interpretable parameters, offering theoretical guarantees and good extrapolation capacity. Data-driven models, in contrast, employ flexible function approximators such as neural networks, excelling in empirical performance but requiring large datasets and often lacking interpretability. Each approach has its limitations: physics-based models struggle when governing equations are unknown or difficult to specify, while data-driven models perform poorly in low-data regimes and fail to extrapolate beyond training distributions.

Recent hybrid approaches attempt to bridge this gap. Physics-Informed Neural Networks (PINNs) \cite{RPK} impose physical constraints as soft penalties in the loss function, while methods like SINDy \cite{BPK} discover governing equations from data. Neural operators \cite{Lu,NO} learn mappings between function spaces. These approaches, however, remain fundamentally constrained: they attempt to satisfy both physical laws and data fidelity within a single model, leading to competing optimization objectives that can be difficult to balance effectively.

We introduce Hybrid-Cooperative Learning (HYCO), a fundamentally different approach that trains two models in parallel: a physical model grounded in differential equations (with unknown parameters $\Lambda$) and a synthetic model guided by data (with parameters $\Theta$). Rather than encoding physics as constraints within the synthetic model, HYCO couples the models through an interaction loss that measures the discrepancy between their predictions. The key mechanism can be understood through three loss components. First, the physical model has its own loss $\L_{phy}(\Lambda)$, measuring how well it fits observational data (if any is provided). Second, the synthetic model has a loss $\L_{syn}(\Theta)$ quantifying its fit to the data. Third, and most importantly, the interaction loss $\L_{int}(\Theta, \Lambda)$ measures the distance between the two models' predictions across the domain. The overall objective combines these terms as
\[
\min_{\Theta, \Lambda}\,  \alpha \L_{syn}(\Theta) + \beta \L_{phy}(\Lambda) + \L_{int}(\Theta, \Lambda) ,
\]
where $\alpha, \beta \geq 0$ control the relative importance of each component. Crucially, either weight may be zero, allowing scenarios where one model trains exclusively through the interaction term without direct access to observational data.

As a simple diagnostic baseline, one may set the interaction term to zero, which decouples the two components and reduces HYCO to two independently trained models.

This cooperative formulation admits a natural game-theoretic interpretation. Instead of jointly optimizing both models, we can view them as two autonomous agents or players. The physical model seeks to minimize its own cost \begin{equation}
\label{l1}
\L_1 := \beta\L_{phy}(\Lambda) + \L_{int}(\Theta, \Lambda),
\end{equation}
while the synthetic model aims to minimize 
\begin{equation}
\label{l2}
\L_2 := \alpha\L_{syn}(\Theta) + \L_{int}(\Theta,\Lambda).
\end{equation}
Each player strives to reduce its own loss while being constrained by the other's choices through the shared interaction term. The solution concept is a Nash equilibrium: a configuration where neither agent can unilaterally reduce its loss by changing parameters. This perspective naturally suggests an alternating minimization scheme: at each iteration, we update the physical parameters $\Lambda$ while holding $\Theta$ fixed, then update $\Theta$ while holding $\Lambda$ fixed. This back-and-forth process nudges both models toward mutual agreement, enabling them to cooperate despite having potentially different objectives.

\smallskip
\noindent\textbf{Contributions.}
This short paper has a primarily expository and methodological scope. Our main contributions are:
\begin{itemize}[leftmargin=*,itemsep=2pt]
\item[$\diamond$] \emph{Symmetric hybridization.} We formulate HYCO as a genuinely \emph{peer-to-peer} coupling of a physics-based model and a data-driven model, avoiding the usual paradigm where physics is only enforced as a constraint on the synthetic component.
\item[$\diamond$] \emph{Game-theoretic viewpoint and practical algorithm.} We reinterpret the coupled training as a two-player game with alternating updates, leading to a simple and implementable training procedure.
\item[$\diamond$] \emph{Stochastic interaction via ghost points.} We introduce a lightweight interaction-evaluation strategy based on randomly sampled ``ghost points'', which decouples the interaction loss from the physical mesh and is convenient in sparse/irregular-data settings.
\end{itemize}

\subsection{Advantages of HYCO} The scope of HYCO is broad. At its simplest, it can be used to embed standard machine learning models with physical biases, much like in the spirit of PINNs. However, the symmetric treatment of both models in HYCO provides notable flexibility: neither component is subordinate to the other, and the framework imposes minimal restrictions on model choice. 

Besides, \emph{parallelization} emerges naturally from HYCO's design. During each training iteration, both models can compute their predictions and gradients simultaneously, with synchronization required only when evaluating the interaction loss. This parallel structure enables distributed training across multiple devices, and in principle could offer computational advantages, though the need to train two models simultaneously means efficiency gains depend significantly on implementation details and the specific problem structure.

\emph{Privacy awareness} is another practical advantage. The two components need not exchange raw data or model parameters; in our implementations they can interact through predicted states (and/or residual-type summaries). This is particularly attractive in federated or multi-institution settings where data cannot leave its source due to regulatory or proprietary constraints. A full privacy analysis (e.g., potential information leakage through shared predictions) is beyond the scope of this short note.

Finally, HYCO is designed to be robust to data sparsity. The interaction loss provides implicit regularization: even when observational data is limited, the physical model benefits from the synthetic model's learned patterns, while the synthetic model gains stability from physical consistency. This mutual regularization can be particularly valuable for extrapolation beyond the training domain, a scenario where purely data-driven approaches often struggle.

\section{Methodology}
\label{sec:methodology}

\subsection{Problem Formulation}

We consider a general setting where observational data is available from some underlying ground truth model or physical system. The goal is to simultaneously recover the solution field and identify unknown parameters characterizing the system. To this end, we employ two models operating in parallel.

The \emph{physical model} is based on a structural ansatz with unknown parameters $\Lambda$. This could be a discretized differential equation solver, a system of ordinary differential equations, or any model grounded in domain knowledge where certain coefficients, forcing terms, or boundary conditions are unknown. Given parameters $\Lambda$, the physical model produces a prediction $\u_{phy}$. If observational data is available to the physical model, it seeks to minimize
\begin{equation}
\label{eq:loss_phy}
\L_{phy}(\Lambda) := \frac{1}{M} \sum_{i=1}^M \ell(\u_{phy}(x_i), \u^D(x_i)) + \mathscr{P}(\Lambda),
\end{equation}
where $\u^D$ represents observational data at locations $\{x_i\}_{i=1}^M$, $\ell$ is a loss function (typically mean squared error), and $\mathscr{P}(\Lambda)$ is an optional regularization term on the parameters.

The \emph{synthetic model} is a data-driven component, typically a neural network with parameters $\Theta$. It learns to approximate the underlying system directly from data by minimizing
\begin{equation}
\label{eq:loss_syn}
\L_{syn}(\Theta) := \frac{1}{M} \sum_{i=1}^M \ell(\u_{syn}(x_i), \u^D(x_i)) + \mathscr{P}(\Theta).
\end{equation}
The key innovation lies in the \emph{interaction loss}, which couples the two models by measuring the discrepancy between their predictions over the entire domain $\Omega$:
\begin{equation}
\label{eq:loss_int}
\L_{int}(\Theta, \Lambda) := \int_{\Omega} \|\u_{syn}(x) - \u_{phy}(x)\|^2 dx.
\end{equation}
This term enforces consistency between the models across the full spatial domain, not merely at the observed data points. The overall optimization problem becomes
\begin{equation}
\label{eq:main_objective}
\min_{\Theta, \Lambda} \, \alpha \L_{syn}(\Theta) + \beta \L_{phy}(\Lambda) + \L_{int}(\Theta, \Lambda),
\end{equation}
where $\alpha, \beta \geq 0$ are weights that balance the objectives. Importantly, either weight may be zero, allowing one model to train exclusively through the interaction term without direct access to observational data.

The framework is agnostic to whether the problem is static or dynamic; the essential structure remains unchanged. In particular, this formulation extends naturally to time-dependent problems, where the domain becomes $\Omega \times [0,T]$ and all integrals are taken over both space and time. 

\subsection{The HYCO Algorithm}

Rather than jointly minimizing \eqref{eq:main_objective}, we adopt an alternating optimization scheme motivated by the game-theoretic perspective. Recalling $\L_1, \L_2$ defined at \eqref{l1}-\eqref{l2}, at each iteration we update one set of parameters while holding the other fixed, taking a gradient descent step to minimize either the physical or the synthetic loss:

\medskip
\begin{algorithm}[H]
\caption{HYCO Alternating Optimization}
\label{alg:hyco}
\SetAlgoLined
\KwIn{Initial parameters $\Lambda^{(0)}, \Theta^{(0)}$; weights $\alpha, \beta$; learning rates $\eta_\Lambda, \eta_\Theta$; maximum iterations $K$}
\KwOut{$(\Lambda^{(K)}, \Theta^{(K)})$}
\For{$k = 0, 1, 2, \ldots, K-1$}{
    \textbf{Physical update:} 
    $\Lambda^{(k+1)} \leftarrow \Lambda^{(k)} - \eta_\Lambda \nabla_\Lambda \L_1$\;
    
    \textbf{Synthetic update:} 
    $\Theta^{(k+1)} \leftarrow \Theta^{(k)} - \eta_\Theta \nabla_\Theta \L_2$\;
    
    \If{stopping criterion satisfied}{
        \textbf{break}\;
    }
}
\end{algorithm}

\subsection{Ghost Points and Efficient Computation}

A critical implementation detail concerns the computation of the interaction loss $\L_{int}$. Evaluating the integral in \eqref{eq:loss_int} exactly at every training iteration would require discretizing the entire domain $\Omega$ on a fine grid, which can be prohibitively expensive, especially for high-dimensional problems or when the physical model itself already operates on a dense mesh.

To address this, we introduce the concept of \emph{ghost points}. The term `ghost" reflects that these points are not part of the dataset and do not directly contribute to the data-fitting losses $\L_{phy}$ or $\L_{syn}$. Rather, they are auxiliary evaluation points sampled randomly from $\Omega$ at each training iteration, serving solely to approximate the interaction loss. At iteration $k$, we sample $H$ points $\{x_h\}_{h=1}^H$ uniformly (or according to some other sampling strategy) from $\Omega$ and approximate
\[
\L_{int}(\Theta, \Lambda) \approx \frac{|\Omega|}{H} \sum_{h=1}^H \|\u_{syn}(x_h) - \u_{phy}(x_h)\|^2.
\]
This stochastic approximation offers several advantages. First, it decouples the computational cost of evaluating $\L_{int}$ from the discretization resolution of the physical model. Second, the random sampling introduces a form of stochastic regularization that can help prevent overfitting to specific regions of the domain. Third, different ghost points are sampled at each iteration, ensuring that over the course of training, the interaction loss is enforced approximately uniformly across $\Omega$. This approach draws inspiration from Random Batch Methods and stochastic optimization techniques commonly used in large-scale machine learning.

The number of ghost points $H$ represents a trade-off: too few points may yield noisy gradient estimates, while too many increase computational cost without commensurate benefit. In our experiments, we found that moderate values of $H$ (on the order of hundreds) suffice to obtain stable training.

\subsection{Stopping Criterion}

A natural question is when to terminate the alternating optimization. One possible approach monitors the change in the interaction loss between successive iterations. If the models have reached approximate agreement, the interaction loss should stabilize. We can thus define a stopping criterion as
\begin{equation}
\label{eq:stopping}
\frac{|\L_{int}^{(k+1)} - \L_{int}^{(k)}|}{\L_{int}^{(k)}} < \epsilon,
\end{equation}
where $\epsilon > 0$ is a tolerance parameter. When this condition is satisfied, further alternating updates are unlikely to significantly improve the solution, and training can be halted.

We emphasize that this is only one possible implementation choice. Other stopping criteria could be based on the change in parameter values, validation loss, or a fixed computational budget. The appropriate choice depends on the specific application and available computational resources.

\section{Related Work}
\label{sec:related}

\noindent
Hybrid modeling approaches have emerged as a promising direction to combine physics and data. Physics-Informed Neural Networks (PINNs) \cite{RPK} encode PDE residuals directly in the loss function, enabling the network to learn solutions that approximately satisfy governing equations. Variants include variational formulations \cite{KZK} and Bayesian extensions \cite{BPINNs}. Neural operators \cite{Lu,NO} learn mappings between function spaces, offering a different paradigm for data-driven PDE solving. SINDy \cite{BPK} discovers governing equations from data through sparse regression.

These methods share a common limitation: they attempt to satisfy both physical constraints and data fidelity within a single model, leading to challenging multi-objective optimization. In contrast, HYCO decouples these objectives by employing two independent models that coordinate through consistency constraints. This design is conceptually closer to ensemble learning \cite{ZM}, though aggregation occurs during training rather than at inference. The alternating optimization bears similarity to ADMM-based PINN methods \cite{YXH}, but HYCO does not impose a hierarchical relationship between models.

HYCO's philosophy also connects to federated learning \cite{MMRHA}, where distributed agents train collaboratively while preserving privacy. In HYCO, the physical and synthetic models exchange only predictions, not parameters or raw data, making the framework naturally privacy-aware. Recent work on minimal variance model aggregation \cite{BO} and bilevel optimization for neural operators \cite{ZXL} explores related ideas, though from different perspectives.

A distinguishing feature of HYCO is the symmetric treatment of the two models. Neither is subordinate to the other; instead, they cooperate as peers to reach consensus. This symmetry provides flexibility in model selection and allows for scenarios in which one model is trained solely through interaction (e.g., by setting $\beta = 0$ or $\alpha = 0$ in the minimization problem), a configuration that is not naturally accommodated by most hybrid frameworks. For comparison, we mention the works \cite{PAK, XU}, where hybrid models for inverse problems are developed, but with a clear hierarchy between the physical and the synthetic model.

\section{Numerical Experiments}
\label{sec:experiments}

\noindent
We validate HYCO on two benchmark problems that are well-established in the inverse problems literature: the Gray-Scott reaction-diffusion system \cite{GrayScott} and the heterogeneous Helmholtz equation. The former tests parameter identification in a dynamic setting with complex pattern formation, whereas the latter tests spatial parameter reconstruction with limited observational data.

\subsection{Gray-Scott Reaction-Diffusion System}

The Gray-Scott model \cite{GrayScott} describes autocatalytic chemical reactions between two species $u = u(t,x,y)$ and $v = v(t,x,y)$ on the plane, and it is well known to exhibit rich spatio-temporal dynamics and pattern formation. On the unit square $\Omega = [0,1]^2 \subset \R^2$ with periodic boundary conditions, the system of PDEs reads
\begin{equation}
\label{eq:gray_scott}
\begin{dcases}
u_t - D_u \Delta u - uv^2 + F(1 - u) = 0, \\
v_t - D_v \Delta v + uv^2 - (F + k)v = 0, \\
u(0,x,y) = u_0(x,y), \quad v(0,x,y) = v_0(x,y).
\end{dcases}
\end{equation}
Here, $D_u$ and $D_v$ are the diffusion coefficients for $u$ and $v$, respectively; $F$ is the feed rate of $u$; and $k$ denotes the decay rate of $v$. The ground-truth parameters used to generate the data in this experiment are:
\[
D_u = 2\times10^{-6}, \quad D_v = 0.8\times10^{-6}, \quad F = 0.018, \quad k = 0.051.
\]
The initial values are
\[
u_0(x,y) = \chi_{B_{0.1}(1/2)}, \quad
v_0(x,y) = \chi_{\Omega \setminus B_{0.1}(1/2)},
\]
where $\chi_A$ denotes the characteristic function of the set $A$ and $B_{0.1}(1/2)$ is the ball of radius 0.1 centered at $(x,y) = (1/2,1/2)$. In what follows, we assume the initial datum is known. This information will be supplied to all compared models.

The dataset is constructed by numerically solving the Gray-Scott equations on a 64$\times$64 grid, using a finite difference scheme with small temporal steps ($5000$ time steps from $0$ to $T=2000$) to generate a reference true solution $\u^D$. From this solution, we sample the values of both species at $M = 5000$ random points, obtaining the dataset
\[
\D := \{\u^D(x_i,y_i,t_i)\}_{i=1}^{M} = \{u^D(x_i,y_i,t_i),v^D(x_i,y_i,t_i)\}_{i=1}^{M},
\]
where the superscript $D$ denotes ``Data". Note, in particular, that our data is not equispaced in space nor in time.

The goal of this experiment is to recover the time-evolution of the solution over the entire domain $\Omega$, generalizing to points not included in the training set, and to recover the diffusivity parameters $D_u, D_v$, henceforth assumed unknown.

\smallskip
\noindent
\textbf{Model configurations.} We compare three methods: HYCO, a baseline deep neural network, and a physics-informed neural network (PINN). The HYCO framework consists of two cooperatively trained components:
\begin{itemize}[leftmargin=*,itemsep=2pt]
\item[$\diamond$] \textit{Physical model} ($\u_{phy}$): This consists of a finite difference discretization of \eqref{eq:gray_scott} on the same 64$\times$64 spatial grid. The unknown diffusion coefficients $\Lambda = \{D_u, D_v\}$ are initialized at $D_u = 1\times 10^{-6}$ and $D_v = 0.5\times10^{-6}$, deliberately chosen away from the true values to test parameter recovery.
\item[$\diamond$] \textit{Synthetic model} ($\u_{syn}$): This is a feedforward neural network with four hidden layers (128 neurons each, ReLU activation), mapping $(x,y,t) \in \R^3$ to $(u_{syn}, v_{syn}) \in \R^2$. This network learns to approximate the solution directly from data while being regularized by the physical model through the interaction loss.
\end{itemize}
The HYCO interaction loss is computed using $H=1000$ ghost points at each iteration.
The baseline neural network uses the same architecture as HYCO's synthetic component but is trained purely on data without physical constraints. The PINN also uses the same architecture but with $\tanh$ activation (required for computing PDE residuals via automatic differentiation). It is trained with 50,000 randomly sampled collocation points to enforce the PDE residuals in the loss function, and also attempts to identify $D_u$ and $D_v$ by treating them as additional learnable parameters. HYCO and the baseline network were trained for 600 epochs. The PINN required 3500 epochs to achieve convergence.

\begin{figure}[h!]
    \centering
    \includegraphics[width=0.7\linewidth]{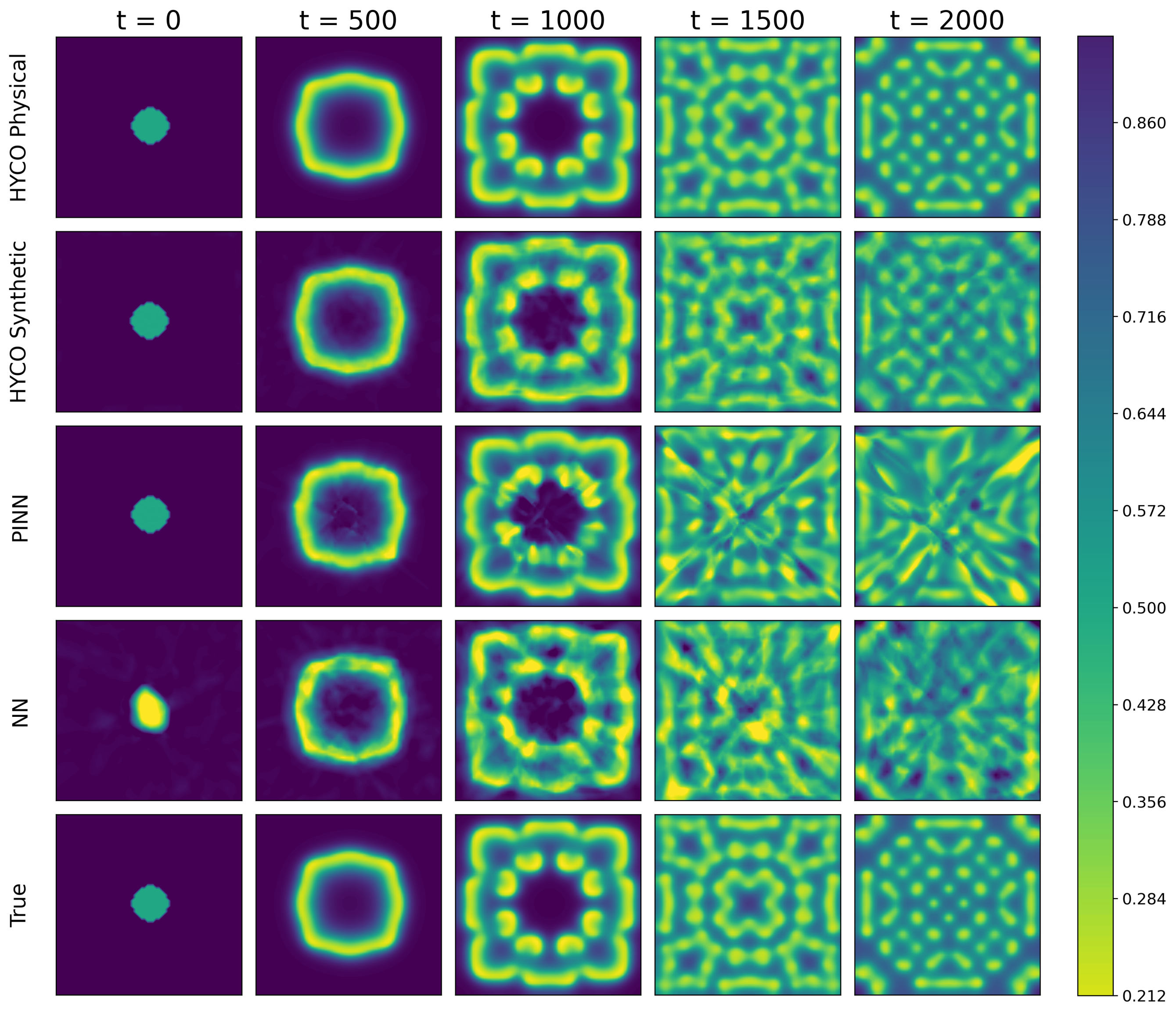}
    \caption{Gray-Scott experiment: evolution of the $u$-component. Rows show HYCO Physical, HYCO Synthetic, pure NN, PINN, and ground truth.}
    \label{fig:GS}
\end{figure}

\smallskip
\noindent
\textbf{Results.} Figure~\ref{fig:GS} shows the evolution of the $u$-component at five time snapshots. The physical and synthetic components of HYCO both capture the pattern formation dynamics throughout the time horizon. The baseline neural network and PINN initially approximate the solution but become unstable as pattern complexity increases, failing to extrapolate beyond the temporal range of the training data.

Table~\ref{tab:gs} reports the final normalized $L^2$ errors 
\[\mathsf{e}^m_s := \frac{\|\u-\u_m\|_{L^2(\Omega\times[0,T])}}{\|\u\|_{L^2(\Omega\times[0,T])}},\]
computed on the full discretization grid. Here $\u$ is the ground truth solution and $\u_m$ the prediction of the model $m$ (HYCO physical or synthetic, PINN or NN). The HYCO physical model achieves error 0.020, followed by HYCO synthetic (0.063). The PINN and baseline network exhibit higher errors (0.123 and 0.153 respectively) due to temporal instability.

\begin{table}[h!]
\centering
\renewcommand{\arraystretch}{1.2}
\begin{tabular}{|c|c|c|c|}
\hline
\textbf{HYCO Physical} & \textbf{HYCO Synthetic} & \textbf{PINN} & \textbf{NN} \\
\hline
\textbf{0.020} & 0.063 & 0.123 & 0.153 \\
\hline
\end{tabular}
\caption{Normalized $L^2$ errors for Gray-Scott approximations.}
\label{tab:gs}
\end{table}

For parameter identification, Table~\ref{tab:gsparams} compares the recovered diffusion coefficients. HYCO accurately identifies the ground-truth values, while the PINN converges to incorrect parameters despite fitting the training data reasonably well. This occurs because the PINN's flexible network can compensate for incorrect parameters by learning a solution that approximately satisfies the PDE with biased coefficients. The HYCO framework mitigates this issue through the consistency constraint between the physical and synthetic models.

\begin{table}[h!]
\centering
\renewcommand{\arraystretch}{1.2}
\begin{tabular}{|c|c|c|}
\hline
 & \textbf{$D_u$} & \textbf{$D_v$} \\
\hline
\textbf{Ground truth} & $2\times10^{-6}$ & $0.8\times10^{-6}$ \\
\hline
\textbf{HYCO Physical} & $\mathbf{1.989\times10^{-6}}$ &$\mathbf{0.799\times10^{-6}}$ \\
\hline
\textbf{PINN} & $0.0999\times10^{-6}$ & $0.0381\times10^{-6}$ \\
\hline
\end{tabular}
\caption{Identified diffusivity parameters for Gray-Scott.}
\label{tab:gsparams}
\end{table}


\subsection{Helmholtz Equation: Static Inverse Problem}

The heterogeneous Helmholtz equation is a classical benchmark for inverse problems in the machine learning and PDE literature, widely used for testing both forward and inverse solvers. Our setup closely follows Example 3.3.2 from the work by Yang, Meng, and Karniadakis on Bayesian PINNs \cite{BPINNs}.

The governing equation is the 2D heterogeneous Helmholtz equation with null Dirichlet boundary conditions, namely
\begin{equation}
\label{eq:helmholtz}
\begin{dcases}
- \nabla \cdot (\kappa(x, y) \nabla \u) + \eta(x, y)^2 \u = f(x, y), & \text{for } (x, y) \in \Omega, \\
\u(x, y) = 0, & \text{for } (x, y) \in \partial\Omega,
\end{dcases}
\end{equation}
where $\Omega = [-\pi, \pi]^2$. Here $\u = \u(x,y)$ is the unknown solution field, $\kappa(x,y)$ represents the diffusion coefficient, $\eta(x,y)$ is the wave number, and $f(x,y)$ is a known forcing term. The ground-truth coefficients are defined as
\begin{equation*}
\kappa(x, y) = \varphi(x,y;\alpha_1,c_1) + 1, \quad \eta(x, y) = \varphi(x,y;\alpha_2,c_2) + 1,
\end{equation*}
with parameters
\[
\alpha_1 = 4, \quad \alpha_2 = 1, \quad c_1 = (-1, -1), \quad c_2 = (2, 1),
\]
and $\varphi(\cdot,\cdot;\alpha,c)$ being the Gaussian function with amplitude $\alpha$ and center $c = (c_x, c_y)$:
\[
\varphi(x,y;\alpha,c) = \alpha e^{-(x-c_x)^2 - (y - c_y)^2}.
\]
Following \cite{BPINNs}, the reference solution is set to
\[
\u(x,y) = \sin(x)\sin(y),
\]
with the corresponding forcing term $f(x,y)$ computed by substituting this solution into \eqref{eq:helmholtz}.

The dataset is constructed by randomly placing $M=25$ sensors within specified regions of the domain and evaluating the reference solution at these locations:
\[
\D := \{\u^D(x_i,y_i)\}_{i=1}^{M}.
\] 
The forcing term $f$ is assumed known and supplied to all models.

The goal of this experiment is to solve the inverse problem of recovering the parameters $\Lambda := \{ \alpha_1, c_1, \alpha_2, c_2\} \in \R^6$ that define the coefficients $\kappa$ and $\eta$. Since $c_1$ and $c_2$ are 2D points, the total number of parameters to identify is six.

\smallskip
\noindent
\textbf{Data generation.} We test three scenarios that progressively reduce the observational coverage:
\begin{itemize}[leftmargin=*,itemsep=2pt]
\item[$\diamond$] \textit{Full domain ($\Omega$)}: sensors randomly placed across the entire domain $[-\pi, \pi]^2$.
\item[$\diamond$] \textit{Partial domain ($\mathsf{Q}_1$)}: sensors restricted to $[-\pi/2, \pi]^2$ (approximately three-quarters of the domain).
\item[$\diamond$] \textit{Quarter domain ($\mathsf{Q}_2$)}: sensors restricted to $[0, \pi]^2$ (one-quarter of the domain, most challenging).
\end{itemize}
In each case, $M=25$ sensors are randomly positioned within the specified region.

\smallskip
\noindent
\textbf{Model configurations.} We compare three methods: HYCO, a classical finite element inverse solver (FEM), and a physics-informed neural network (PINN). The HYCO framework consists of two cooperatively trained components:
\begin{itemize}[leftmargin=*,itemsep=2pt]
\item[$\diamond$] \textit{Physical model} ($\u_{phy}$): This consists of a finite element discretization of \eqref{eq:helmholtz} with unknown parameters $\Lambda = \{\alpha_1, c_1, \alpha_2, c_2\}$ initialized randomly. In this experiment the physical model is also trained on the dataset $\D$, minimizing the data mismatch loss along with the interaction term.
\item[$\diamond$] \textit{Synthetic model} ($\u_{syn}$): This is a feedforward neural network with ResNet-like architecture featuring two hidden layers (256 neurons each, ReLU activation) and skip connections, mapping $(x,y) \in \R^2$ to $\u_{syn}(x,y) \in \R$. 
\end{itemize}
The HYCO interaction loss is computed using $H=200$ ghost points at each iteration.

The finite element method (FEM) baseline solves the classical inverse problem by optimizing the six parameters via gradient-based methods on the discretized PDE, minimizing only the data mismatch. In the partial-coverage settings, this inverse problem becomes markedly ill-posed: many parameter configurations can explain the measurements inside the observed region while producing very different fields elsewhere, so the optimization is sensitive to initialization and regularization. The PINN uses a similar neural architecture (with $\tanh$ activation required for automatic differentiation), enforcing the PDE residual through 400 randomly sampled collocation points (boundary and interior), and treating the six parameters as learnable alongside the network weights. All models were trained for a maximum of 2000 epochs using the Adam optimizer.

\begin{table}[h!]
\centering
\renewcommand{\arraystretch}{1.2}
\begin{tabular}{|c|c|c|c|}
\hline
\textbf{Domain} & \textbf{Model} & \textbf{$\mathsf{e}^m_p$} & \textbf{$\mathsf{e}^m_s$} \\
\hline
\multirow{3}{*}{$\Omega$} 
& FEM            & 0.024  & 0.012 \\
& HYCO Physical  & \textbf{0.016} & \textbf{0.012} \\
& PINN           & 0.158  & 0.061 \\
\hline
\multirow{3}{*}{$\mathsf{Q}_1$}
& FEM            & \textbf{0.024}   & \textbf{0.012} \\
& HYCO Physical  & 0.031  & \textbf{0.012} \\
& PINN           & 0.178   & 0.078 \\
\hline
\multirow{3}{*}{$\mathsf{Q}_2$}
& FEM            & 1.109    & 0.366 \\
& HYCO Physical  & \textbf{0.041}   & \textbf{0.072} \\
& PINN           & 0.750    & 0.361 \\
\hline
\end{tabular}
\caption{Helmholtz results: parameter error $\mathsf{e}_p$ and solution error $\mathsf{e}_s$ across observational domains.}
\label{tab:helm}
\end{table}

\noindent
\textbf{Results.} Table~\ref{tab:helm} summarizes the parameter error $\mathsf{e}^m_p$ and the solution error $\mathsf{e}^m_s$ across the three test cases, where we report $\mathsf{e}^m_p:=\|\Lambda^{(m)}-\Lambda^{\dagger}\|_2/\|\Lambda^{\dagger}\|_2$ and $\mathsf{e}^m_s:=\|u^{(m)}-u^{\dagger}\|_{L^2(\Omega)}/\|u^{\dagger}\|_{L^2(\Omega)}$ (approximated on the reference grid). When data covers the full domain ($\Omega$), all methods achieve accurate parameter recovery with small errors. As observational coverage decreases, however, performance differences become pronounced. For the intermediate case ($\mathsf{Q}_1$), HYCO and FEM perform comparably, while the PINN exhibits larger errors. For the most challenging case ($\mathsf{Q}_2$), HYCO achieves parameter error $\mathsf{e}^m_p = 0.041$ and solution error $\mathsf{e}^m_s = 0.072$, substantially outperforming both FEM ($\mathsf{e}^m_p = 1.109$) and PINN ($\mathsf{e}_p = 0.750$).

\begin{figure}[h!]
    \centering
    \includegraphics[width=0.7\linewidth]{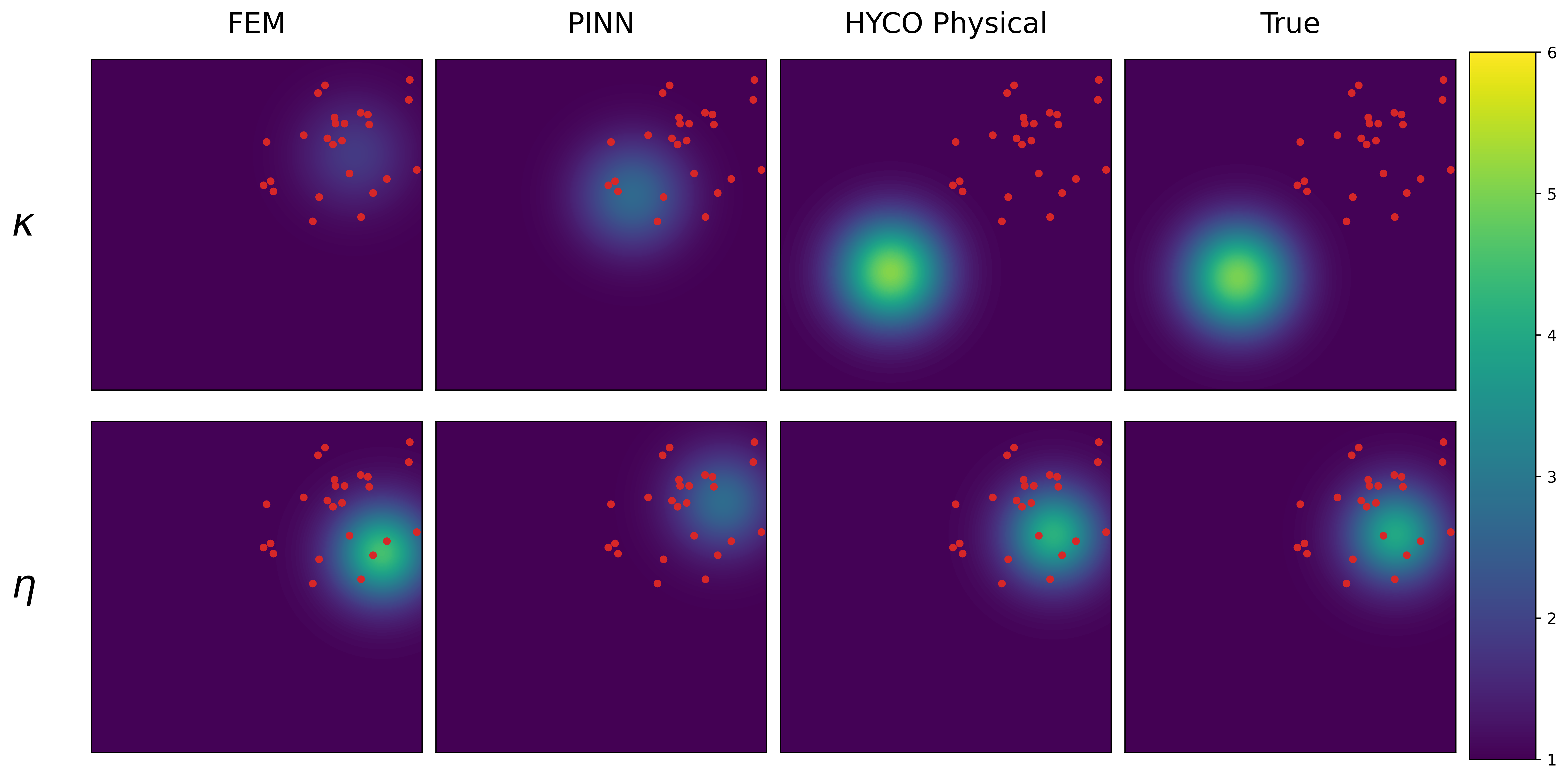}
    \caption{Recovered coefficients for Helmholtz equation with observations in $\mathsf{Q}_2$: $\kappa$ (top row), $\eta$ (middle row). Columns show FEM, PINN, HYCO Physical, and ground truth. Red dots indicate sensor positions.}
    \label{fig:helm_param}
\end{figure}

Figure~\ref{fig:helm_param} visualizes the recovered coefficients for the most challenging case ($\mathsf{Q}_2$). HYCO accurately reconstructs both $\kappa$ and $\eta$ across the entire domain, including regions without observations (marked by red dots). In contrast, FEM and PINN converge to alternative parameter configurations that fit the observed data well but fail to generalize outside the observation region. This demonstrates that all three models find mathematically valid solutions to the inverse problem in the sense that they reproduce the observed data. However, HYCO's cooperative framework, which implicitly encourages consistency between the physical and synthetic models beyond the data region, leads to parameter estimates closer to the ground truth.

\begin{figure}[t!]
    \centering
    \includegraphics[width=0.7\linewidth]{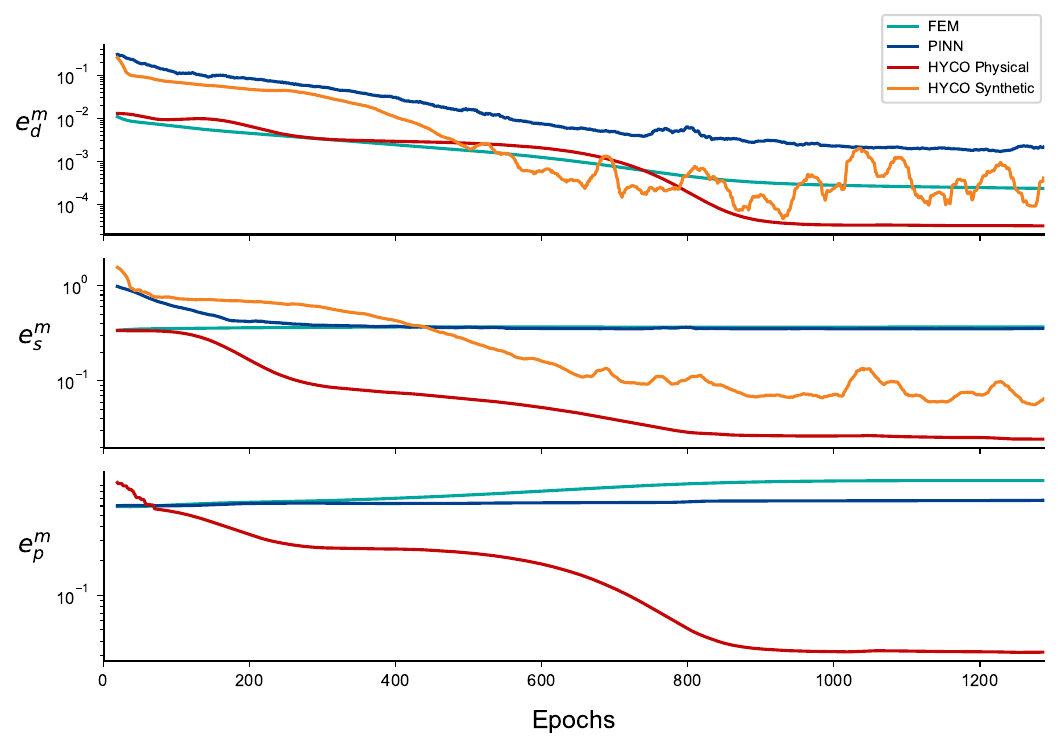}
    \caption{Error evolution for Helmholtz experiment with observations in $\mathsf{Q}_2$. Top: data mismatch. Middle: solution error. Bottom: parameter error.}
    \label{fig:helm_error}
\end{figure}

Figure~\ref{fig:helm_error} shows the evolution of errors during training for the $\mathsf{Q}_2$ case. All models effectively reduce the data mismatch (top panel), which is expected since this term is explicitly minimized. However, only HYCO maintains low solution and parameter errors throughout training (middle and bottom panels), while FEM and PINN converge to solutions with accurate data fit but incorrect parameters.

\section{Conclusions and Future Work}
\label{sec:conclusions}

\noindent
This paper presents HYCO (Hybrid-Cooperative Learning), a modeling strategy that integrates physical and synthetic models through joint optimization while encouraging alignment between their predictions. The work presented here is based on the extended preprint \cite{XXX}, where a comprehensive theoretical and numerical analysis is developed.

In HYCO, both models are treated as active participants in the learning process, with mutual influence serving as a form of regularization. Through numerical experiments we demonstrated that HYCO can perform well on tasks ranging from solution recovery to parameter identification. These advantages are particularly evident when observations are sparse or spatially localized, where the interplay between models helps reconstruct missing information.

However, important limitations and open questions remain. The computational cost of training two models simultaneously means that efficiency gains from parallelization depend significantly on implementation details and problem structure. While our experiments show competitive performance, the overhead of maintaining two models is non-negligible. Furthermore, the theoretical analysis remains incomplete: convergence guarantees for the game-theoretical formulation are established in \cite{XXX} only for simplified convex settings, and the behavior in general non-convex cases (which dominate practical applications) requires further investigation.

Looking ahead, several research directions emerge. HYCO's general architecture invites application to more complex scenarios: multiphysics problems where each model represents a different physical process, domain decomposition where models govern different spatial regions, and real-world applications in climate modeling, robotics, or medical imaging. Testing HYCO with advanced architectures such as transformers or neural operators as synthetic models remains an important validation step. On the theoretical front, extending convergence analysis to non-convex settings and investigating potential privacy vulnerabilities (as highlighted in federated learning literature) are critical challenges. Finally, while the game-theoretic perspective offers interpretability, further exploration of this connection may yield improved training algorithms or relaxation gap results analogous to recent work in neural network optimization.

In summary, HYCO represents an ongoing research effort that treats hybrid modeling through the lens of cooperative agents rather than constrained optimization. While promising results have been demonstrated on benchmark problems, substantial work remains to fully understand the method's capabilities, limitations, and optimal deployment strategies.

\section*{Acknoledgements}

\noindent
This work was supported by SURE-AI: The Norwegian Centre for Sustainable, Risk-Averse, and Ethical AI, grant 357482, Research Council of Norway. E.Z. was supported by the Grant PID2023-146872OB-I00-DyCMaMod of MICIU (Spain), as well as by the ERC Advanced Grant CoDeFeL(ERC-2022-ADG-101096251).

\bibliographystyle{plain}
\bibliography{bib}

\end{document}